\documentclass{amsart}
\usepackage[T1]{fontenc}
\usepackage{amssymb}
\usepackage{amssymb,booktabs}
\usepackage[hyperref]{degt}

\def\2{\color{red}}
\def\3{\color{green}}
\def\4{\color{blue}}

\def\Dg:{{\bf Dg:\enspace}\ignorespaces}

\def\KK-{\pdfstr{K3}{$K3$}-\penalty0\ignorespaces}

\def\0{}

\def\F{\Bbb F}
\def\FF{\bar{\F}}

\def\spn{\operatorname{span}}
\def\V{\operatorname{V}}

\def\cl[#1]{[\![#1]\!]}

\def\NS{\operatorname{NS}}

\def\resp.{resp\PERIOD}

\let\leq\le

\let\M=M
\let\o=o

\let\ex=e
\let\iso=p

\def\CC{\Cal C}

\def\i{^\mathrm{i}}
\def\ii{^\mathrm{ii}}

\newcommand{\UU}{\hat{u}}
\newcommand{\VV}{\hat{v}}
\newcommand{\krL}{\mathcal L}

\newcounter{row}
\def\therow{\arabic{row}}

\def\vtab{\vbox\bgroup\def\(##1){[##1]}%
 \ialign\bgroup\hss\refstepcounter{row}\ \therow##\hss\quad&
  \strut\hss$##$\hss\quad&$\singset{##}$\ \hss\cr
  \noalign{\hrule\vspace{2pt}}
  \omit\ \hss$\#$\hss\quad&\text{ord}&\omit\text{pattern}\hss\cr
  \noalign{\vspace{1pt}\hrule\vspace{2pt}}
  }
\def\endvtab{\crcr\noalign{\vspace{1pt}\hrule}\egroup\egroup}

\pdef\two{\root3\of2}

\title[$K3$ complete intersections with maximal number of lines]
{On $K3$ complete intersections with maximal number of lines }

\author{Alex Degtyarev}

\address{%
Department of Mathematics\\
Bilkent University\\
06800 Ankara, TURKEY}

\email{
degt@fen.bilkent.edu.tr}

\thanks{%
A.D. was partially supported by the T\"{U}B\DOTaccent{I}TAK grant 123F111%
}

\author{S\l awomir Rams}
\address{Faculty of Mathematics and Computer Science, Jagiellonian University,
	ul.\ {\L}ojasiewicza 6,  30-348 Krak\'{o}w, Poland}
\email{slawomir.rams@uj.edu.pl}

\thanks{%
    S.R. was supported by the National Science Centre, Poland, Opus  grant
	no.\ 2024/\allowbreak53/\allowbreak B/\allowbreak ST1/\allowbreak01413
}



\begin{document}
\newcommand{\smooth}{}
\begin{abstract}
We give an explicit description of
$K3$-sextics and $K3$-octics with the maximal number of lines.  As a consequence, we show that these numbers ($42$ and $36$, respectively) are also attained in any characteristic $p>3$. In characteristic $5$, we find a supersingular octic with $38$ lines.
\end{abstract}

\maketitle

\section{Introduction}\label{S.intro}

Numerous results on
configurations  of lines
(or higher degree rational curves)
 on various classes of surfaces
are spread over the literature with many open questions.
It is due to the progress of the last decade
that our understanding
 of
configurations of lines on polarized complex \KK-surfaces became almost complete
(see \cite{segre, rams.schuett, DIS, degt:lines,degt:sextics,
Degtyarev.Rams.triangular, Degtyarev.Rams, Degtyarev.Rams.sextics, shioda56,
Miyaoka:lines, rams.schuett:quasi} and further references therein). In
particular, for  complex \KK-surfaces that are complete intersections in
a projective space we know that
\begin{itemize}
	\item Schur's quartic is the only degree-$4$ surface that contains $64$ lines,
	\item
$K3$-sextics with $42$ lines constitute an irreducible $1$-parameter
family,
	\item
there are exactly two $K3$-octics that contain $36$ lines,
	\end{itemize}
and no $K3$-surface of degree $4$ (resp.\ $6$, resp.\ $8$)
contains more than $64$ (resp.\ $42$, resp.\ $36$) lines.
(Here and below, all uniqueness statements are modulo the appropriate
projective linear group.)

We know
a great deal on the large configurations of lines
on smooth degree-$4$ $K3$-surfaces over
other  fields.
However, in the case of positive characteristic our techniques are
limited by the complexity of Torelli-type results for $K3$-surfaces.
On the other hand, one can use  \cite{deligne-lift, NO} (see also \cite{JangLift})  and the results on complex $K3$-surfaces
	to study the positive characteristic case (see
\cite{rams.schuett:char2, rams.schuett:char3, degt:linessupersin}
for an illustration of this approach for quartics).
That is why the
knowledge of
explicit equations of complex $K3$ surfaces  with many lines plays
important r\^{o}le in the study of the positive characteristic case.

The main subject of this note are generators of ideals of complex  $K3$ complete intersections with maximal number of lines, but since the surfaces we consider are defined over $\QQ$,
we 
discuss line configurations on their reductions modulo $p$ as well
(see Remarks~\ref{rem-sextic-pos-char},~\ref{rem-octic-prime-pos-char} and
\ref{rem-octic-bis-pos-char}).

The equation of Schur's quartic~\cite{Schur:quartics} goes back to the XIXth
century,
whereas those of the remaining extremal $K3$ complete intersections
remain unknown.  Here we fill this gap
and
prove the following theorems.


\theorem[see \autoref{proof.sextics}] \label{thm.sextics}
Let $Q \subset  \PP^4(\C)$ be the singular quadric
$\V(x_1 x_2 - x_3 x_4)$
and
assume that
$\lambda^3\ne0,4$.
\roster
\item\label{sextics.def}
The
$K3$-sextic $X_{\lambda} \subset \PP^4(\C)$ given by the intersection
of $Q$ with the cubic
\begin{equation} \label{eq-sextic42lines}
	x_0^3 + x_1^3 + x_2^3 + x_3^3 + x_4^3
	{{}-3\lambda}x_0 x_1 x_2 = 0
\end{equation}
is smooth.  Each of the two
	families of $2$-planes
	on $Q$
 endows
$X_{\lambda}$ with an elliptic fibration
{\rm(}by planar cubics{\rm)}
of
type
$6 \tA_2 \oplus 3 \tA_1$.
In particular, $X_{\lambda}$ contains $42$ lines.
\item\label{sextics.iso}
Two sextics $X_\lambda$ and $X_\mu$ are projectively equivalent if and
	only if $\lambda^3=\mu^3$. In particular,
	the pencil $\{X_{\lambda}\}$ is not isotrivial.
\item\label{sextics.all}
Hence, any $K3$-sextic with $42$ lines is projectively equivalent to one of
	$X_\lambda$.
\endroster
\endtheorem

If $\lambda^3 \in \{0, 4\}$,  the cubic \eqref{eq-sextic42lines} meets the quadric $Q$
	along a sextic with simple singularities.
Such a surface  contains at most $36$
lines by \cite[Theorem~1.1]{Degtyarev.Rams.sextics}.


For $K3$-octics with $36$ lines we have
fewer possibilities by  \cite[Theorem~1.2]{degt:lines} (see also
\cite[Theorem~1.1]{Degtyarev.Rams}).

\theorem[see \autoref{proof.octics}] \label{thm.octics}
Every complex $K3$-octic in $\PP^5(\C)$ that contains $36$ lines is
projectively equivalent
either to the octic
$X':=\V(q'_0, q'_1, q'_2)$ given by
\begin{align} \label{eq-quadrcs'}
	q'_0&=x_0^2+x_1^2+2x_2x_3-x_5^2, \nonumber \\
	q'_1&=x_2^2+x_3^2+2x_0x_1+x_4^2,\\
	q'_2&=2x_0x_2-x_4x_5. \nonumber
\end{align}
or to the octic $X'':=\V(q''_0, q''_1, q''_2)$
given by
\begin{align}  \label{eq-quadrcsb}
	q''_0&=2x_0x_3+x_1x_2-x_1x_5-x_2x_4+3x_4x_5, \nonumber \\
	q''_1&=2x_1x_4+x_2x_0-x_2x_3-x_0x_5+3x_5x_3,\\
	q''_2&=2x_2x_5+x_0x_1-x_0x_4-x_1x_3+3x_3x_4. \nonumber
\end{align}
\rom(In full agreement with \cite{degt:lines,Degtyarev.Rams}, extremal octics
are triquadrics.\rom)
\endtheorem

Apart from the obvious beauty of
these equations
(\eg, one cannot overlook the similarity of the construction of $\{X_{\lambda}\}$ with
the Barth quintic),
one should observe that their reductions
give us
examples of $K3$-sextics and octics with many lines over fields of characteristic $p>3$.
In  \autoref{discrim-1458}  we also present explicit models of complex $K3$-sextics with minimal discriminants
(\cf. \cite[Theorem~1.5.3]{degt:smoothmodels}).
Therefore, apart from
an attempt to complete the picture of the geometry of varieties that were studied already by Klein
(\cf. \cite[\S\,10.3.3]{dolgachev:cag}), our interest in explicit
equations
of $K3$ complete intersections with largest configurations of lines
is motivated by the fact that
these equations can
be useful for further research.

The paper is organized
as follows. In
\autoref{proof.sextics}
 we prove \autoref{thm.sextics} on sextics,
discuss the N\'{e}ron--Severi groups of some members of the pencil $\{X_{\lambda}\}$ (see \autoref{discrim-1458})
and their reductions modulo $p \neq 3$ (see \autoref{rem-sextic-pos-char}).
Octics are treated in \autoref{S.octics}:
we prove \autoref{thm.octics} and
discuss the reductions modulo $p$
of complex $K3$-octics with $36$ lines (see Remarks~\ref{rem-octic-prime-pos-char},~\ref{rem-octic-bis-pos-char}).
In particular, we construct a smooth $K3$-octics with $38$ lines over $\FF_5$.

\convention\label{conv}
Unless otherwise stated,  in this note we work over the field $\C$ of complex
numbers, and we abbreviate $\Cp{n}:=\Cp{n}(\C)$.
Lines are
degree-one curves in the ambient
projective space.
We reserve the notation $\zeta_m:=\exp(2\pi i/m)$.
\endconvention



\newcommand{\PPi}{\Pi}
\newcommand{\PPiA}{\Pi'}

\newcommand{\PPiB}{\Pi''}
\newcommand{\FA}{F'}
\newcommand{\FB}{F''}
\newcommand{\RA}{R'}
\newcommand{\RB}{R''}

\newcommand{\ddiscr}{\det}

\section{ $K3$-sextics with $42$ lines} \label{proof.sextics}

We maintain the notation of {\autoref{S.intro}}
and  assume that
$\lambda^3\ne0,4$.
Consider
$\PP^4(\C)$ with homogenous coordinates $(x_0: \ldots :x_4)$ and the point
$P_0 := (1:0:\ldots:0)$.

Recall that
the
quadric $Q=\{x_1x_2-x_3x_4=0\}\subset\PP_4$ is  the cone with vertex
$P_0$ over
a smooth quadric $S=\Cp1\times\Cp1\subset\Cp3$: the latter is given by
the same equation in coordinates $(x_1:\ldots:x_4)$ for $\Cp3$.
Obviously, the rulings
\[*
\RA_t=\{x_3=tx_1,\ x_2=tx_4\},\quad \RB_t=\{x_4=tx_1,\ x_2=tx_3\}, \quad t \in \C,
\]
of $S$
endow
$Q$
with two families of $2$-planes
$$
\PPiA_t=\operatorname{span}(\{P_0\} \cup \RA_t),\quad
\PPiB_t=\operatorname{span}(\{P_0\} \cup \RB_t), \quad t \in \C.
$$
Then, each smooth $K3$-sextic $X := Q\cap\V(f)$, where $f$ is a
degree-$3$ polynomial
carries the  elliptic fibrations
{$|F^*_t|:=\{\Pi^*_t \cap X \,|\, t \in \PP_1\}$, ${*}=\prime$ or
$\prime\prime$.}

\proof[Proof of  \autoref{thm.sextics}]
To begin with, observe that, up to projective equivalence
{(\eg, $x_0\mapsto{\pm\zeta_3} x_0$),
$X_\lambda$
does depend on $\lambda^3$ only; that is why it is $\lambda^3$ that appears in
most equations.

Consider the plane $\PPiA_t$ with
coordinates 
  $(x_0:u:tv:tu:v)$.
The
fiber 
$\FA_t$
is given by
\[ \label{eq-hesse-disguise}
x_0^3
{{}-3\lambda}t\,x_0uv+(1+t^3)(u^3+v^3)=0 .
\]
If $t^3+1=0$, this
planar cubic splits into a line and a conic.
For $t^3 + 1\neq 0$,
pick a branch $\tau:=\root3\of{1+t^3}$ and put
\[*
\UU := \tau u,\quad \VV := \tau v,\quad
\mu:=\lambda t/\tau^2,
\]
so that \eqref{eq-hesse-disguise} becomes
the Hesse cubic
$$x_0^3+\UU^3+\VV^3
{{}-3\mu} x_0 \UU \VV =0.$$
Hence, apart from the three fibers ($t^3=-1$) of
type  $\tA_1$,
each consisting of a
line and
irreducible conic $x_0^2=3\lambda {\4 t} uv$, we obtain six fibers of type $\tA_2$ (triangles of lines)
corresponding to the values
\[*
\mu^3=1,\quad\text{\ie},\quad t^6+(2-\lambda^3)t^3+1=0.
\]
(Recall that we assume $\lambda^3\neq 0,4$.)
This gives us exactly
$21$ lines as components of the singular fibers of $|\FA_t|$. The symmetry
$x_3\leftrightarrow x_4$
shows that $21$ more lines are components of
the nine singular fibers of $|\FB_t|$.

We have already described $24=(21\cdot\text{lines}+3\cdot\text{conics})$
	fiber components of $|\FA_t|$. Any exceptional divisor on $X_\lambda$ would
	have to be an extra fiber component of  this pencil, which is impossible.
	Hence, $X_\lambda$
	is smooth.

To prove statement~\iref{sextics.iso}, we observe that a projective
equivalence would preserve the two pencils above. The singular fibers
$\Pi^*_t\cap Q$ are those with $t^3=-1$ (type~$\tA_1$) or $t^3=\Ga_1,\Ga_2$,
where $\Ga_1,\Ga_2$ are the roots of $\Ga^2+(2-\Gl^3)\Ga+1=0$. The position
of the type~$\tA_1$ fibers fixes the coordinates in the base of the pencil up
to the subgroup $\SG3$ of the M\"{o}bius group generated by
$t\mapsto{\zeta_3} t$ and
$t\mapsto t\1$. They preserve the pair $\Ga_1$,
$\Ga_2=\Ga_1\1$ and, hence, $\lambda^3$.

Finally, statement~\iref{sextics.all} follows from~\iref{sextics.iso} and the
fact~\cite{degt:lines,Degtyarev.Rams} that, up to $\PGL(5,\C)$, smooth
$K3$-sextics with $42$ lines constitute an irreducible $1$-parameter family.
\endproof

By \cite{degt:lines} (see also \cite[Theorem~1.1]{Degtyarev.Rams.sextics}) the configuration of
lines on $X_{\lambda}$ is independent of $\lambda$ (its graph is called $\Psi_{42}$ in \cite{degt:lines}).
Moreover, the sublattice of $\NS(X_{\lambda})$   spanned by the classes of the $42$
lines is  of rank $19$.

 In particular,  \autoref{thm.sextics}\iref{sextics.iso} and \cite{oguiso} imply that $\rho(X_{\lambda}) = 19$ for uncountably many values of the parameter $\lambda$. One can prove  that $$\rho(X_{-8}) = 19.$$
	Indeed, one   checks that $p=19$ and $p=31$ are primes of good reduction for $X_{-8}$ and all
lines on $X_{-8}$
	are $\F_p$-rational for those primes;
then, one follows the  standard approach of \cite{vanluijk} to disprove
	the fact that $X_{-8}$ is singular. Since the description of this technique can
be easily found in the literature,  we leave the details to the reader.

Recall that, by \cite[Theorem~1.5(3)]{degt:smoothmodels}, (modulo a projective equivalence)
there exist exactly six  complex singular $K3$-sextics with $\ls|\ddiscr(\NS)| \leq 48$:
two of which satisfy  $\ddiscr(\NS) = -39$, whereas for the remaining  four we have
$\ddiscr(\NS) = -48$.
Below, we find
three of these surfaces {(those with $42$ lines)} in the pencil $\{X_\lambda\}$.

\remark \label{discrim-1458}
(i) Consider the conic 
$$C= \V\bigl(\two\,x_0-(x_1+x_2+x_3+x_4),\ \ x_1-x_2-i\,(x_3-x_4),\ \ x_1x_2-x_3x_4\bigr).$$
By a direct computation $C \subset{X_{3\two}}$ and  the discriminant of the lattice spanned by the $42$ lines and the conic $C$ is  $(-48)$.
 Thus   \cite[Theorem~1.5(3)]{degt:smoothmodels} yields
$$
\ddiscr(\NS(X_{3\two})) = -48 \, .
$$

\noindent
(ii) Let $(\lambda_0)^3=-\bigl(23\pm7\sqrt{13}\bigr)/2$. A standard
Groebner basis computation with \cite{singular} shows that the sextic $X_{\lambda}$ contains $60$ conics.
Then, one can check that the  discriminant of the lattice spanned by the classes of the $42$ lines and the $60$ conics equals $(-39)$.  This implies
$$
\ddiscr(\NS(X_{\lambda_0})) = -39 \, .
$$
By \autoref{thm.sextics}\iref{sextics.iso}, the two surfaces are not
projectively equivalent. Hence, we do obtain both determinant-minimizing
sextics found in \cite{degt:smoothmodels}.
\endremark

\remark
As
explained above,
our recent paper \cite{Degtyarev.Rams.sextics} shed
new
light on the geometry of $K3$-sextics with $42$ lines. This in turn allowed us to find the
generator \eqref{eq-sextic42lines} via finite field search with help of Gr\"{o}bner bases and lifting its results to characteristic zero. Since the final result of our search is defined over $\Q$, it was sufficient to apply CRT in the process of lifting.
\endremark

Finally, let us apply \autoref{thm.sextics} to obtain  positive-characteristic examples of
$K3$-sextics with $42$ lines.

	\remark \label{rem-sextic-pos-char}
	Let
	$J_6 \subset \ZZ[\lambda, x_0, \ldots, x_4]$
	be the ideal generated by the quadric $(x_1x_2-x_3x_4)$, the cubic polynomial $f_{\lambda}$ defined in \eqref{eq-sextic42lines}
	and all $(2 \times 2)$-minors of the Jacobian matrix of the pair $(x_1x_2-x_3x_4, f_{\lambda})$ with respect to the homogenous variables $x_0$, $\ldots$, $x_4$. Then  a strong Gr\"{o}ebner basis computation (see \cite[Chapter~4]{AdamsL1994}) with \cite{singular}  proves that
	$$
	3 \lambda^2 (\lambda^3 - 4) x_j^7 \in J_6 \mbox{ for } j = 0, \ldots, 4 \, .
	$$
	In particular, this implies that, for a prime
	$p \neq 3$ and every  parameter $\lambda$ such that
	$$
	\lambda (\lambda^3 - 4) \neq 0 \bmod p \, ,
	$$
	the $K3$-sextic $(X_{\lambda})_{\FF_p}$ is smooth.
	Given a prime $p >3$, each smooth sextic $(X_{\lambda})_{\FF_p}$ contains at least $42$
lines (by the standard specialization argument).
On the other hand (see the proof of \autoref{thm.sextics}), every line on
	 $(X_{\lambda})_{\FF_p}$ is a
fiber component of one of
     the elliptic fibrations $|\FA_t|$, $|\FB_t|$.
Thus,  $(X_{\lambda})_{\FF_p}$ contains exactly $42$ lines.

     Similarly, standard analysis of the geometry of  the elliptic fibrations $|\FA_t|$, $|\FB_t|$  shows that every  smooth  $(X_{\lambda})_{\FF_2}$ contains exactly $42$ lines.
	\endremark

\section{ $K3$-octics with $36$ lines (c.f.  \cite[Table~1]{Degtyarev.Rams})} \label{S.octics}

By \cite{degt:lines}
or \cite[Theorem~1.1]{Degtyarev.Rams}, there
are
exactly two $K3$-octics with $36$ lines. Their transcendental lattices and Fano graphs
are found in  \cite[Table~1]{Degtyarev.Rams};
see also \cite[\S\,2.3]{Degtyarev.Rams} for the definition of the Fano
graph of a surface.

\subsection{The $K3$-octic $X'$ with $36$ lines and $\ls|\Aut_hX'|=16$.} \label{ss-octic16}
\newcommand{\OO}{\mathfrak{O}}

We consider the group of automorphisms of $\PP^5(\C)$ generated by the maps
\begin{eqnarray*}
\varphi_1\:(x_0,x_1,x_2,x_3, x_4, x_5) &\mapsto&  (x_0,x_1,-x_2,-x_3, -x_4, -x_5) \, , \\
\varphi_2\:(x_0,x_1,x_2,x_3, x_4, x_5) &\mapsto&  (x_0,x_1,x_2,x_3, i x_5, -i x_4) \, , \\
\psi_{}\:(x_0,x_1,x_2,x_3, x_4, x_5)   &\mapsto&    (x_0,-x_1,i x_2,-i x_3, i x_4, x_5)
\end{eqnarray*}
It is immediate that
$\varphi_1^2 = \varphi_2^2 = \psi^4= (\varphi_2 \psi)^2=\id$
and $\varphi_1$ commutes with
both $\varphi_2$ and~$\psi$.
In other words,
\begin{equation} \label{eq-autXp}
\<\varphi_1, \varphi_2, \psi  \> \simeq
\Z/2\times\DG8.
\end{equation}
Furthermore, the web $\spn(q'_0, q'_1, q'_2)$ with the base locus~$X'$,
see~\eqref{eq-quadrcs'}, is invariant under $\varphi_1, \varphi_2, \psi$.
Consider the points
\[*
\gathered
\alignedat2
P\i_4  &:= \bigl(1:i:0:0:1-i:0\bigr), &\qquad Q\i_4 &:=  \bigl(0 : 0 : 1 : i : 0 : 1+i) , \\
P\i_8  &:= \bigl(1:i:0: -i \sqrt{2+2i}: \sqrt{2}: 0\bigr), & Q\i_8  &:= \bigl(0: \sqrt{2+2i}: 1; -i : 0 : \sqrt{2} )  , \\
P\ii_8 &:= \bigl(1:-i:0: i \sqrt{2-2i}: \sqrt{2}: 0\bigr), & Q\ii_8 &:= \bigl(0: \sqrt{2-2i}: 1; i : 0 : \sqrt{2} )  ,
\endalignedat\\
\aligned
	P\i_{16} &:= \bigl(1:i (1+\sqrt{2})^{1/2}:0:(i-1)(1+\sqrt{2}):0:i(2+2\sqrt{2})^{1/2}\bigr) , \\
	Q\i_{16} &:= \bigl(0:1-i:i(1+\sqrt{2}):(\sqrt{2}-1)^{1/2}:\sqrt{2}:0\bigr)
\endaligned
\endgathered
\]
and define the lines $\ell^i_k := \<P^i_k, Q^i_k \>$. One can easily see
that
$$
\ell^i_k\subset X'\quad\text{for all four applicable pairs $(i,k)$}.
$$
Define $\OO^i_k$
as the orbit of $\ell^i_k$
under the action of
\eqref{eq-autXp}.
It is
straightforward
that
\begin{equation} \label{eq-orbit-length}
\ls|\OO^i_k|=k.
\end{equation}
Moreover,
letting $\krL'$ denote the lattice spanned by the $36$ lines in the four
orbits,
\begin{equation} \label{eq-discrXp}
\rank(\krL') = 20,\qquad
\ddiscr{\krL'} = -32 \, .
\end{equation}

\lemma \label{lem-Xp}
The $K3$-octic $X'=\V(q'_0, q'_1, q'_2)$
as in \eqref{eq-quadrcs'}
contains exactly $36$ lines. Moreover, we have
$$
\Aut_h(X')
=\<\varphi_1, \varphi_2, \psi  \> \simeq\Z/2\times\DG8,\qquad
\operatorname{T}_{X'} =
\begin{bmatrix}
	4 & 0  \\
	0 & 8 \end{bmatrix}.
$$
The set of lines on $X'$ consists of the four orbits of action of
$\Aut_h(X')$, and
the classes of the $36$ lines generate $\NS(X')$.
\endlemma
\proof
By \cite[Theorem~1.2]{degt:lines}
or \cite[Theorem~1.1]{Degtyarev.Rams}, the
$K3$-octic $X'$ contains at most $36$ lines.
Thus, \eqref{eq-orbit-length} implies that ${X':={}}\V(q'_0, q'_1, q'_2)$  is one of the two $K3$-octics listed in  \cite[Table~1]{Degtyarev.Rams}.
Then \eqref{eq-discrXp} implies that
$\det \operatorname{T}_{X'}\mathrel|32$ and, therefore, $X'$ is the octic with the Fano
graph $\Theta'_{36}$ in  \cite[Table~1]{Degtyarev.Rams}.
\endproof

For the sake of completeness let us briefly discuss the reductions  $X'_{\FF_p}$ for various
primes $p$.

	\remark \label{rem-octic-prime-pos-char}
	Let
	$J'_8 \subset \ZZ[x_0, \ldots, x_4]$
	be the ideal generated by the quadrics $q'_0, q'_1, q'_2$,
	and all $(3 \times 3)$-minors of the Jacobian matrix of the triple  $(q'_0, q'_1, q'_2)$ with respect to the homogenous variables $x_0$, $\ldots$, $x_5$.
	We compute
	a strong Gr\"{o}bner basis of $J'_8$ and use the ideal membership test to show that
	$$
	64 \, x_j^6 \in J'_8 \mbox{ for } j = 0, \ldots, 5 \, .
	$$
	Therefore the $K3$-octic $(X')_{\FF_p}$
is smooth for every
	odd prime $p$ and, by  the standard specialization argument
	it contains
	at least $36$ lines.  Finally, one computes  a strong Gr\"{o}bner basis of the ideal of the Fano scheme of the octic $X'$ over $\ZZ$ (more precisely:  strong Gr\"{o}bner bases of ideals of subschemes of the Fano scheme supported on various Schubert cells of the Grassmanian $\operatorname{Gr}(2,6)$)
	to check that there are no extra lines on  $(X')_{\FF_p}$  for odd primes $p \neq 5$. For
	$p=5$ one checks directly that
	the (supersingular) $K3$-octic  $(X')_{\FF_5}$ contains exactly $38$ lines.
	
	The reader, who prefers more conceptual approach, can  also exploit the arithmetic
 techniques in the spirit of \cite{shimada-supersingular} to check for which primes
$p$ the octic $(X')_{\FF_p}$ is supersingular
	and can acquire extra lines.
The fact that the rank-$20$ lattice $\NS(X')$ is  generated by the classes of lines
	simplifies some arguments.
	\endremark

\subsection{The $K3$-octic $X''$ with $36$ lines and $\ls|\Aut_hX''|=144$.}

Consider the group of automorphisms of $\PP^5(\C)$ generated by
\begin{eqnarray*}
	\sigma\:(x_0, x_1, x_2, x_3, x_4, x_5) &\mapsto&  (x_1,x_2,x_0,x_4, x_5, x_3) \, , \\
	  \tau\:(x_0, x_1, x_2, x_3, x_4, x_5) &\mapsto&  (x_0, \zeta_3 x_1, \zeta_3^2 x_2, x_3, \zeta_3 x_4, \zeta_3^2 x_5)\,.
\end{eqnarray*}
Then $\sigma^3 = \tau^3 =
\operatorname{id}$,
$\sigma \tau = \tau \sigma$ so that
\begin{equation} \label{eq-heisenberg}
\langle \tau, \sigma \rangle \simeq (\Z/3)^2 \, .
\end{equation}
Furthermore, the web $\spn(q''_0, q''_1, q''_2)$ with the base locus~$X''$,
see~\eqref{eq-quadrcsb}, is invariant under the action of~\eqref{eq-heisenberg}.
Consider the points}
\[*
\gathered
\alignedat2
P_1 &:= \bigl(1:-2:1:0:0:1\bigr),& \qquad Q_1 &:= \bigl(0:0:3:1:1:-2\bigr) \, , \\
P_2 &:= \bigl(1:1:-2:0:1:0\bigr),&        Q_2 &:= \bigl(0:0:3:1:-2:1\bigr)  \, ,
\endalignedat\\
\aligned
	P_{3} &:= \bigl(2: -1 - \sqrt{3}: -1 + \sqrt{3}:0:1 - \sqrt{3}:1 + \sqrt{3}\bigr)  \, , \\
	Q_{3} &:= \bigl(0: 3 + 3\sqrt{3}: 3 - 3\sqrt{3}: 2: -1 + \sqrt{3}: -1 - \sqrt{3}\bigr)  \, , \\
	P_{4} &:= \bigl(2: -1 + \sqrt{3}: -1 - \sqrt{3}:0:1 + \sqrt{3}:1 - \sqrt{3}\bigr)  \, , \\
	Q_{4} &:= \bigl(0: 3 - 3\sqrt{3}: 3 + 3\sqrt{3}: 2: -1 - \sqrt{3}: -1 + \sqrt{3}\bigr)  \,
\endaligned
\endgathered
\]
and define the four lines
$\ell_k := \<P_k, Q_k\>$.
A direct computation shows that, for each $k=1,\ldots,4$, the line $\ell_k$
lies on~$X''$ and the orbit of $\ell_k$ under the action
of~\eqref{eq-heisenberg} consists of exactly nine lines
$$
\ell_{q,r,k} := \sigma^q \tau^r(\ell_k), \quad  q,r = 0, 1, 2.
$$
This time, for the lattice $\krL''{}\subset\NS(X'')$
spanned by the classes of the $36$ lines
 we have  $\rank(\krL'') = 19$.
However,
 \begin{equation} \label{eq-discrXb}
 \ddiscr\bigl(\spn(\krL'', h)\bigr) = -36,  \quad \text{where $h$ is the hyperplane section}.
 \end{equation}


\lemma \label{lem-Xb}
The $K3$-octic $X''=\V(q''_0, q''_1, q''_2)$,
see \eqref{eq-quadrcsb},
contains exactly $36$ lines. Moreover, the set of lines on $X''$ consists of
the four orbits of the action of~\eqref{eq-heisenberg}
and $\spn(\krL'', h) = \NS(X'')$. Finally, we have
$$
\operatorname{T}_{X''} =
\begin{bmatrix}
	6 & 0  \\
	0 & 6 \end{bmatrix}
$$
\endlemma

\proof
Similar to \autoref{lem-Xp},
the statement is a direct consequence of~\eqref{eq-discrXb} and
\cite[Table~1]{Degtyarev.Rams}. The configuration of lines is $\Theta''_{36}$
in the notation of \latin{loc.\ cit}.
\endproof


	\remark \label{rem-octic-bis-pos-char}
	Let
	$J''_8 \subset \ZZ[x_0, \ldots, x_4]$
	be the ideal generated by the quadrics $q''_0, q''_1, q''_2$,
	and all $(3 \times3)$-minors of the Jacobian matrix of the triple  $(q''_0, q''_1, q''_2)$ with respect to the homogenous variables $x_0$, $\ldots$, $x_5$.
	We use
	a strong Gr\"{o}bner basis of $J'_8$
and use the ideal membership test to show that
	$$
	216 \, x_j^6 \in J''_8 \mbox{ for } j = 0, \ldots, 5 \, .
	$$
	This shows that the $K3$-octic $(X'')_{\FF_p}$ is smooth for
every
prime $p>3$, and
	it contains
	at least $36$ lines. Then  a strong Gr\"{o}bner basis computation of the ideal of the Fano
scheme of the octic $X''$ proves that there are no extra lines on  $(X'')_{\FF_p}$  for
primes $p > 3$.
	\endremark

\remark
In contrast to the sextics in \autoref{thm.sextics}, the search for octics with
$36$ lines is far more complicated.
(To begin with, the results of \cite{Degtyarev.Rams.sextics} dictate that
the quadric~$Q$ in \autoref{thm.sextics} must have a single singular point.
Hence, it can be assumed fixed, whereupon one is essentially working with
hypersurfaces rather than complete intersections.)
We followed the strategy designed in \cite{veniani:equations}
(building on  the classical results of Nikulin \cite{Nikulin:forms}), \ie,
we studied the automorphisms of the N\'{e}ron--Severi and transcendental lattices
of $X'$ and $X''$
described in \cite{degt:lines,Degtyarev.Rams}
to arrive at constraints on the quadrics  that contain the octics in
question,
which \emph{are} triquadrics by \latin{loc.\ cit}.
We omit
technical details.
\endremark

\remark
Let us mainatin the notation of \eqref{eq-heisenberg} and introduce two  additional
automorphisms of $\PP^5(\C)$:
\begin{eqnarray*}
	\iota\:(x_0, x_1, x_2, x_3, x_4, x_5) &\mapsto&    (x_0, x_2, x_1, x_3, x_5, x_4) \, , \\
	\theta\:(x_0, x_1, x_2, x_3, x_4, x_5) &\mapsto& (3 x_3, 3 x_4, 3 x_5, x_0, x_1, x_2) \,
\end{eqnarray*}
Then $\iota^2=\theta^2=(\iota \sigma)^2=(\iota \tau)^2=
\operatorname{id}$, and $\theta$ commutes with the three automorphisms  $\sigma$, $\tau$, $\iota$.
One can easily see that
\begin{equation}\label{eq-18group}
	\langle \tau, \sigma, \iota, \theta  \rangle \simeq
	\bigl((\Z/3)^2\rtimes\Z/2\bigr)\times\Z/2  \subset \operatorname{Aut}_h(X'').
\end{equation}
According to
\cite[Theorem~1.2]{degt:lines} we have $\ls|\operatorname{Aut}_h(X'')|=144$.
(In \cite[Theorem~1.1]{Degtyarev.Rams} it is erroneously stated
that $\ls|\operatorname{Aut}_h(X'')|=72$.)
Found in \eqref{eq-18group} are but $36$ automorphisms. For the sake of completeness let us 
define $t:=1+(\sqrt3-1)\zeta_3$ and 
consider the automorphisms of $\PP^5(\C)$ defined by the following matrices:
$$
\Phi=\begin{pmatrix}
	\sqrt{-3} & \sqrt{-3} & \sqrt{-3} & 3 & 3 & 3\\
	\sqrt{-3} & \sqrt{-3}\,\zeta_3 & \sqrt{-3}\,\zeta_3^2 & 3 & 3\zeta_3 & 3\zeta_3^2\\
	\sqrt{-3} & \sqrt{-3}\,\zeta_3^2 & \sqrt{-3}\,\zeta_3 & 3 & 3\zeta_3^2 & 3\zeta_3\\
	1 & 1 & 1 & \sqrt{-3} & \sqrt{-3} & \sqrt{-3}\\
	1 & \zeta_3 & \zeta_3^2 & \sqrt{-3} & \sqrt{-3}\,\zeta_3 & \sqrt{-3}\,\zeta_3^2\\
	1 & \zeta_3^2 & \zeta_3 & \sqrt{-3} & \sqrt{-3}\,\zeta_3^2 & \sqrt{-3}\,\zeta_3
\end{pmatrix}  \,  ,
$$
and
$$
\Psi=\begin{pmatrix}
	t & t\zeta_3^2 & t\zeta_3^2 & -3 & -3\zeta_3^2 & -3\zeta_3^2\\
	t\zeta_3 & t\zeta_3 & t\zeta_3^2 & -3\zeta_3 & -3\zeta_3 & -3\zeta_3^2\\
	t\zeta_3 & t\zeta_3^2 & t\zeta_3 & -3\zeta_3 & -3\zeta_3^2 & -3\zeta_3\\
	1 & \zeta_3^2 & \zeta_3^2 & -t & -t\zeta_3^2 & -t\zeta_3^2\\
	\zeta_3 & \zeta_3 & \zeta_3^2 & -t\zeta_3 & -t\zeta_3 & -t\zeta_3^2\\
	\zeta_3 & \zeta_3^2 & \zeta_3 & -t\zeta_3 & -t\zeta_3^2 & -t\zeta_3
\end{pmatrix}
.
$$
One checks that
$$
	\langle \tau, \sigma, \iota, \theta, \Phi, \Psi  \rangle =  \operatorname{Aut}_h(X'')
	\quad \mbox{and} \quad \ls|\operatorname{Aut}_h(X'')|=144.
$$

\endremark

\subsection{Proof of \autoref{thm.octics}}\label{proof.octics}
The
statement is an immediate consequence of
\cite[Theorem~1.2]{degt:lines},  \autoref{lem-Xp}, and
\autoref{lem-Xb}.
\qed

\section{Acknowledgements}
Claude Opus 5.5
(Anthropic) was used for  programming in \cite{singular} and discussions on \cite{AdamsL1994}, \cite{veniani:equations} by S.R, who would like to thank Anthropic
for providing him with free access to its Opus~5.5 model.

{
\let\.\DOTaccent
\def\cprime{$'$}
\bibliographystyle{amsplain}


\begin{thebibliography}{10}
	
\bibitem{AdamsL1994}
	W. W. Adams and P. Loustaunau  \emph{An Introduction to {G}r\"{o}bner Bases},
	Grad. Stud. Math \textbf{3},
American Mathematical Society, Providence, RI, 1994

	
	

	
	\bibitem{Conway.Sloane}
	J.~H. Conway and N.~J.~A. Sloane, \emph{Sphere packings, lattices and groups},
	Grundlehren der Mathematischen Wissenschaften [Fundamental Principles of
	Mathematical Sciences], vol. 290, Springer, 1988, With contributions by E.
	Bannai, J. Leech, S. P. Norton, A. M. Odlyzko, R. A. Parker, L. Queen and B.
	B. Venkov. \MR{920369 (89a:11067)}
	
	
	\bibitem{singular}
	W. Decker, G.-M. Greuel, G. Pfister, H. Sch\"{o}nemann: Singular [4.3.2] — A computer algebra system for polynomial computations. https://www.singular.uni-kl.de (2023)
	
	\bibitem{degt:lines}
	A. Degtyarev, \emph{Lines on {S}mooth {P}olarized {K}3-{S}urfaces}, Discrete
	Comput. Geom. \textbf{62} (2019), no.~3, 601--648. \MR{3996938}
	
		\bibitem{degt:smoothmodels}
	\bysame \emph{Smooth models of singular $K3$-surfaces}.
	Rev. Mat. Iberoam. 35 (2019), no. 1, 125--172.
	
	\bibitem{degt:linessupersin}
	\bysame
	\emph{Lines in supersingular quartics},
	J. Math. Soc. Japan \textbf{74}, No.~3, 973--1019 (2022).
	
	
	
	
	\bibitem{degt:sextics}
	\bysame, \emph{Tritangents to smooth sextic curves}, Ann. Inst. Fourier
	(Grenoble) \textbf{72} (2022), no.~6, 2299--2338.
	
	
	
	
	\bibitem{DIS}
	A. Degtyarev, I. Itenberg, and A.~S. Sert\"{o}z, \emph{Lines on quartic
		surfaces}, Math. Ann. \textbf{368} (2017), no.~1-2, 753--809. \MR{3651588}
	
	\bibitem{Degtyarev.Rams}
	A. Degtyarev and S. Rams, \emph{Counting lines with {V}inberg's
		algorithm}, Rev. Mat. Iberoam. \textbf{42}  (2026), No. 2, 489--544.
	
	\bibitem{Degtyarev.Rams.triangular}
	\bysame, \emph{Lines on {$K3$}-quartics via triangular sets}, Discrete Comput.
	Geom. \textbf{73} (2025), 785--814. \MR{3996938}
	
	
		\bibitem{Degtyarev.Rams.sextics}
	\bysame, \emph{Lines on {$K3$}-sextics with simple singularities}, preprint, arXiv:2512.08018 [math.AG] (2026).


      \bibitem{deligne-weil} P. Deligne \emph{La conjecture de Weil pour les surfaces $K3$},
      Invent. Math. 15 (1972), 206--226.


\bibitem{deligne-lift}
P. Deligne,
\emph{Rel\'{e}vement des surfaces $K3$
en caract\'{e}ristique nulle.} in: "Algebraic surfaces (Orsay, 1976–78)",
Lecture Notes in Math. {\bf 868},
Springer, Berlin, 1981, 58--79.


\bibitem{dolgachev:cag} I. V. Dolgachev,
\emph{Classical algebraic geometry. A modern view.}
Cambridge: Cambridge University Press  (2012).

	
	\bibitem{GAP4}
	The GAP~Group, \emph{{GAP -- Groups, Algorithms, and Programming, Version
			4.15.1}}, 2025.
	
	
	
	\bibitem{JangLift} J. Jang, Neron-Severi group preserving lifting of K3 surfaces and applications, Math. Res. Lett. {\bf 22} (2015), 289--802. 
	
	\bibitem{Miyaoka:lines}
		Y. Miyaoka, \emph{Counting lines and conics on a surface}, Publ. Res. Inst. Math.
	Sci. \textbf{45} (2009), no.~3, 919--923. \MR{2569571}
	
	\bibitem{Nikulin:forms}
	V.~V. Nikulin, \emph{Integer symmetric bilinear forms and some of their
		geometric applications}, Izv. Akad. Nauk SSSR Ser. Mat. \textbf{43} (1979),
	no.~1, 111--177, 238, English translation: Math USSR-Izv. 14 (1979), no. 1,
	103--167 (1980). \MR{525944 (80j:10031)}
	
	\bibitem{oguiso} K. Oguiso, Local families of $K3$ surfaces and applications, J. Algebraic Geom. 12 (2003), no. 3,
	405--433.
	
	
	
\bibitem{NO}  N. O. Nygaard and A. Ogus, Tate conjecture for K3 surfaces of finite
	height. Ann. of Math., {\bf 122} (1985), 461--507.
	
	
	\bibitem{rams.schuett:char3}
		S. Rams and M. Sch\"{u}tt, \emph{112 lines on smooth quartic surfaces (characteristic 3)}, Q. J.
	Math. \textbf{66} (2015), no.~3, 941--951. \MR{3396099}
	
	\bibitem{rams.schuett}
	\bysame, \emph{64 lines on smooth quartic surfaces}, Math. Ann. \textbf{362}
	(2015), no.~1-2, 679--698. \MR{3343894}
	
	
	\bibitem{rams.schuett:char2}
	\bysame, \emph{At most 64 lines on smooth quartic surfaces (characteristic 2)}.
	Nagoya Math. J. \textbf{232} (2018), 76--95 .
	
	
	
	\bibitem{rams.schuett:24}
	\bysame, \emph{24 rational curves on {K3} surfaces}, Commun. Contemp. Math.
	\textbf{25} (2023), no.~6, Paper No. 2250008, 23.
	
	\bibitem{rams.schuett:quasi}
	\bysame, \emph{Low degree rational curves on quasi-polarized
		{{\(\mathrm{K}3\)}} surfaces}, J. Pure Appl. Algebra \textbf{229} (2025),
	no.~2, 21 (English), Id/No 107904.
	
	\bibitem{Schur:quartics}
	F. Schur, \emph{Ueber eine besondre {C}lasse von {F}l\"{a}chen vierter
		{O}rdnung}, Math. Ann. \textbf{20} (1882), no.~2, 254--296. \MR{1510168}

\bibitem{segre} B. Segre, \emph{The maximum number of lines lying on a quartic surface}. Q. J. Math. Oxf. Ser. 14 (1943), 86--96.


\bibitem{shimada-supersingular} I. Shimada, \emph{Transcendental lattices and supersingular reduction lattices
	of a singular $K3$-surface}, Trans. Amer. Math. Soc. \textbf{361} (2009), no.~2, 909--949.




\bibitem{shioda56}
I. Shimada and  T. Shioda,
\emph{On a smooth quartic surface containing $56$ lines which is isomorphic as a $K3$ surface to the Fermat quartic},
Manuscr. Math. \textbf{153}(2017), No. 1-2, 279--297.




	

\bibitem{vanluijk} R. van Luijk,
\emph{K3 surfaces with Picard number one and infinitely many rational points.}
Algebra Number Theory \textbf{1} (2007),  1--17.


\bibitem{veniani:equations} D. C. Veniani,
\emph{Symmetries and equations of smooth quartic surfaces with many lines.} (English) Zbl 1444.14092
Rev. Mat. Iberoam. \textbf{36} (2020), No. 1, 233--256.

	
\end{thebibliography}
}
\end{document}